\documentclass[11pt]{article}
\usepackage[a4paper]{anysize}\marginsize{3.5cm}{3.5cm}{1.3cm}{2cm}
\pdfpagewidth=\paperwidth \pdfpageheight=\paperheight
\usepackage{pgf,tikz,hyperref}
\usetikzlibrary{arrows}
\usepackage{amsfonts,amssymb,amsthm,amsmath,eucal}
\usepackage{graphicx}

\theoremstyle{plain}
\newtheorem{thm}{Theorem}
\newtheorem{theorem}[thm]{Theorem}

\newtheorem{lemma}[thm]{Lemma}

\newtheorem{corollary}[thm]{Corollary}

\theoremstyle{definition}

\newtheorem{problem}[thm]{Problem}

\newtheorem{thevarthm}[thm]{\varthmname}

\newenvironment{varthm*}[1]{\trivlist\item[]{\bf #1.}\it}{\endtrivlist}

\newcommand\eps{\varepsilon}
\newcommand\be{\begin{eqnarray*}}
\newcommand\ee{\end{eqnarray*}}

\newcommand\C{\mathbb C}
\newcommand\call{\mathcal L}

\renewcommand\P{\mathbb P}

\newcommand\newop[2]{\def#1{\mathop{\rm #2}\nolimits}}
\newop\edim{edim}
\newop\Zeroes{Zeroes}
\newop\Jac{Jac}
\newop\Ass{Ass}
\newop\SL{SL}
\newop\PGL{{\P}GL}
\newop\Km{Km}
\newop\mult{mult}

\newcommand\keywords[1]{{\renewcommand\thefootnote{}\footnotetext{\textit{Keywords:} #1.}}}
\newcommand\subclass[1]{{\renewcommand\thefootnote{}\footnotetext{\textit{Mathematics Subject Classification (2010):} #1.}}}

\def\endproof{\hspace*{\fill}\endproofsymbol\endtrivlist}

\def\endproofsymbol{\frame{\rule[0pt]{0pt}{6pt}\rule[0pt]{6pt}{0pt}}}

\begin{document}

\author{Jakub Kabat and Beata Strycharz-Szemberg}
\title{Diminished Fermat-type arrangements\\ and unexpected curves}
\date{\today}
\maketitle
\thispagestyle{empty}

\begin{abstract}
   The purpose of this note is to present and study a new series of the so-called unexpected curves.
   They enjoy a surprising property to the effect that their degree grows to infinity,
   whereas the multiplicity at a general fat point remains constant, equal $3$, which
   is the least possible number appearing as the multiplicity of an unexpected curve
   at its singular point. We show that additionally the BMSS dual curves inherits the
   same pattern of behaviour.
\end{abstract}

\subclass{14C20, 14N10, 14N20}


\section{Introduction}
   The notion of unexpected curves has been introduced by Cook II, Harbourne, Migliore and Nagel in \cite{CHMN}.
   Motivated by an example described first by Di Gennaro, Ilardi, and Vall\'es in \cite{GIV}, they observed
   that there exist configurations $Z$ of points in $\P^2$ such that imposing an additional
   point $P$ of multiplicity $m$ on the linear system of curves of degree $d$ vanishing at all points of $Z$ imposes
   less conditions than expected. If this happens, we say that $Z$ \emph{admits an unexpected curve} of degree $d$
   with a point of multiplicity $m$, or shortly, that $Z$ has the $U(2;d,m)$ property. If the dimension of
   the ambient space is clear to be $2$, in order to alleviate notation, we drop the first index and speak simply
   of the $U(d,m)$ property.

   After the appearance of \cite{CHMN} in 2017, the subject has attracted a lot of attention.
   The foundational article \cite{CHMN} studies sets $Z$ with the $U(d,d-1)$ property. Bauer, Malara, Szemberg, and Szpond
   in \cite{BMSS} discovered the first example of an unexpected surface in $\P^3$, more precisely,
   they identified a set $Z$ with the $U(3;4,3)$ property. In the same paper, they observed a phenomenon
   which is now known as the BMSS duality. Harbourne, Migliore, Nagel and Teitler study in \cite{HMNT}
   more systematically higher dimensional unexpected hypersurfaces, in particular a series of such hypersurfaces attached to root systems. Unexpected cones are studied in \cite{HMTG}. Other constructions of unexpected hypersurfaces appear in \cite{ChiMig, DFHMSTG, FGHM}.
   In a series of papers  \cite{Szp19, Szp19c,Szp19multi}, Szpond shows examples
   of unexpected hypersurfaces with multiple fat points and studies series of examples
   building upon Fermat-type arrangements of hyperplanes. Our study is related to Fermat-type
   arrangements as well.


\section{Initial example}
\label{sec: initial}
   In this section, we analyse in detail an example underlying our construction.
   Let $F$ be the Fermat cubic curve given by the defining equation
   $$x^3+y^3+z^3=0.$$
   It is well-known that $F$ has $9$ inflection points. The coordinates of these points
   can be computed explicitly intersecting $F$ with its Hessian curve $H$, which in this case splits in the union of lines given by the equation
   $$xyz=0.$$
   Thus the inflection points on $F$ have coordinates
   $$A_1=(1:-1:0),\;\; A_2=(1:-\eps:0),\;\ A_3=(1:-\eps^2:0),$$
   $$A_4=(1:0:-1),\;\; A_5=(1:0:-\eps),\;\ A_6=(1:0:-\eps^2),$$
   $$A_7=(0:1:-1),\;\; A_8=(0:1:-\eps),\;\ A_9=(0:1:-\eps^2),$$
   where $\eps$ is a primitive root of the unity of order $3$.
   The lines tangent to $F$ at points $A_1,\ldots,A_9$ are given by equations:
   $$\ell_1:\; x+y=0,\;\; \ell_2:\; x+\eps^2y=0,\;\; \ell_3:\; x+\eps y=0,$$
   $$\ell_4:\; x+z=0,\;\; \ell_5:\; x+\eps^2z=0,\;\; \ell_6:\; x+\eps z=0,$$
   $$\ell_7:\; y+z=0,\;\; \ell_8:\; y+\eps^2z=0,\;\; \ell_9:\; y+\eps z=0.$$
   Note that the product of these lines is the polynomial
   \begin{equation}\label{eq: nine lines product}
   g_3=(x^3+y^3)(y^3+z^3)(z^3+x^3).
   \end{equation}
   Consider the arrangement $\call_3=\left\{\ell_1,\ldots,\ell_9\right\}$. It has $3$ triple points in the coordinate points
   $$X_1=(1:0:0),\;\; X_2=(0:1:0),\;\; X_3=(0:0:1)$$
   and $27$ points where only $2$ configuration lines meet.
   We denote the set of these $27$ points by $Y$. The union of the coordinate points is denoted by $X$.
\begin{lemma}\label{lem: 27 points}
   The points in $Y$ form an almost
   complete intersection. The ideal $I(Y)$ of $Y$ is generated in degree $6$ by the following three polynomials
   $$f_1=x^6-y^6,\;\; f_2=x^6-z^6,\;\; f_3=(x^3+z^3)(y^3+z^3).$$
\end{lemma}
\proof
   The first two polynomials define a complete intersection $W_6$ of $36$ points of the form
   $$(1:\tau^{\alpha}:\tau^{\beta}),$$
   where $\tau$ is a primitive root of unity of order $6$ and $\alpha,\beta=1,\ldots,6$.
   Our $27$ points are the set difference of $W_6$ and the complete intersection $W_3$ defined by
   $$g_1=x^3-y^3,\;\; \mbox{ and }\;\; g_2=x^3-z^3.$$
   Obviously, none of points in $W_3$ belongs to the set of zeroes of $f_3$, whereas
   $f_3$ vanishes at all points of $Y=W_6\setminus W_3$.
\endproof
   We are interested in the set
   \begin{equation}\label{eq: Z}
      Z=Y\cup X,
   \end{equation}
   which associated ideal $I(Z)$ is
   $$I(Z)=I(Y)\cap (xy,xz,yz).$$
\begin{lemma}\label{lem: 30 points}
   The ideal $I(Z)$ is generated in degree $7$ by the following polynomials
   $$h_1=x(y^6-z^6),\;\; h_2=y(x^6-z^6),\;\; h_3=z(x^3+y^3)(y^3+z^3),$$
   $$h_4=z(x^3+y^3)(x^3+z^3),\;\; h_5=y(x^3+z^3)(y^3+z^3),\;\; h_6=x(x^3+z^3)(y^3+z^3).$$
\end{lemma}
\proof
   Let $J$ be the ideal generated by $h_1,\ldots,h_6$. It is easy to see that
   $$J\subset I(Z),$$
   i.e., the polynomials $h_i$ for $i=1,\ldots, 6$ vanish at all points of $Z$.
   For the reverse inclusion, let $F\in I(Z)$ be an arbitrary polynomial.
   Since, in particular, $F\in I(Y)$, there are polynomials $\alpha, \beta, \gamma\in\C[x,y,z]$ such that
   \begin{equation}\label{eq: def of F}
      F=\alpha f_1+\beta f_2+ \gamma f_3.
   \end{equation}
   From $F\in I(X)$ we obtain, evaluating at $X_1$, that $\beta(1:0:0)=0$, hence
   there are polynomials $\beta_1, \beta_2 \in\C[x,y,z]$ such that
   \begin{equation}\label{eq: beta}
      \beta=y\beta_1+z\beta_2.
   \end{equation}
   Analogously, evaluating at $X_2$ and $X_3$, we obtain
   \begin{equation}\label{eq: alpha gamma}
      \alpha=x\alpha_1 + z\alpha_2\;\mbox{ and }\; -\alpha+\gamma-\beta=x\gamma_1+y\gamma_2,
   \end{equation}
   for some polynomials $\alpha_1,\alpha_2,\gamma_1,\gamma_2\in\C[x,y,z]$.
   Substituting \eqref{eq: beta} and \eqref{eq: alpha gamma} to \eqref{eq: def of F}, we obtain
   \begin{align*}
      F &= x\alpha_1(y^6-z^6)+z\alpha_2(y^6-z^6)+y\beta_1(x^6-z^6)+z\beta_2(x^6-z^6)\\
      &\qquad +(x\gamma_1+y\gamma_2+x\alpha_1+z\alpha_2+y\beta_1+z\beta_2)(x^3+z^3)(y^3+z^3)\\
      &= \alpha_1 h_1+\beta_1 h_2+ \gamma_1 h_6+\gamma_2 h_5+\alpha_1 h_6 +\beta_1 h_5\\
      &\qquad +z\alpha_2(y^3+z^3)(y^3-z^3+x^3+z^3)+z\beta_2(x^3+z^3)(x^3-z^3+y^3+z^3)\\
      &= \alpha_1 h_1+\beta_1 h_2+ \gamma_1 h_6+\gamma_2 h_5+\alpha_1 h_6 +\beta_1 h_5 +\alpha_2h_3 +\beta_2 h_4.
   \end{align*}
   Thus $F\in J$ which completes the proof.
\endproof
   Lemma \ref{lem: 30 points} yields the following, immediate, consequence.
\begin{corollary}\label{cor: 30 independent conditions}
   The set $Z$ defined in \eqref{eq: Z} imposes independent conditions of forms of degree $7$.
\end{corollary}
\proof
   Indeed, the space of homogeneous polynomials in $\C[x,y,z]$ of degree $7$ has dimension $\binom{9}{2}=36$.
   There are $30$ points in $Z$, so the expected dimension of the vector space of homogeneous
   polynomials of degree $7$ vanishing along $Z$ is $36-30=6$. This is equal to the actual
   dimension established in Lemma \ref{lem: 30 points}.
\endproof
   The main result in this part is the following statement.
\begin{theorem}\label{thm: deg 7}
   The set $Z$ defined in \eqref{eq: Z} has the $U(7,3)$ property.
\end{theorem}
\proof
   Let $P=(a:b:c)$ be a general point in $\P^2$ and let $\Gamma$ be the curve defined
   by
   \begin{align}\label{eq: gamma}
   \begin{split}
      \gamma=\gamma_P(x:y:z) &= a^2(5c^3-a^3)h_1 + b^2(b^3-5c^3)h_2\\
                      &\qquad + c^2(c^3-5a^3)h_3 + c^2(5b^3-c^3)h_4\\
                      &\qquad + 5b^2(a^3-c^3)h_5 + 5a^2(c^3-b^3)h_6.
   \end{split}
   \end{align}
   Obviously $\gamma\in I(Z)$. Moreover $\mult_P\Gamma\geq 3$. This can be checked
   directly computing partial derivatives of $\gamma$ of order $2$ and checking that
   they all vanish at $P$. By the multiple use of the Euler formula this justifies
   the claim. We omit the simple calculations.
\endproof
\subsection{The BMSS dual of $\Gamma$}
   The idea of the BMSS duality is to consider $\gamma$ as a polynomial in variables $a,b,c$
   with parameters $x,y,z$ viewed as coordinates of a general point $Q=(x:y:z)$ in the
   projective plane with the $a,b,c$ coordinates. Thus, reorganizing terms we have
   \begin{align*}
      \gamma=\gamma_Q(a:b:c) &= x(z^6-y^6)a^5+ y(x^6-z^6)b^5+ z(y^6-x^6)c^5 \\
                     &\qquad + 5y(x^3+z^3)(y^3+z^3)a^3b^2- 5x(x^3+z^3)(y^3+z^3)a^2b^3+ \\
                     &\qquad - 5z(x^3+y^3)(z^3+y^3)a^3c^2+ 5x(x^3+y^3)(z^3+y^3)a^2c^3+ \\
                     &\qquad + 5z(x^3+y^3)(x^3+z^3)b^3c^2- 5y(x^3+y^3)(x^3+z^3)b^2c^3.
   \end{align*}
   We claim that $\mult_Q\Gamma\geq 3$. This is again easy to check verifying vanishing
   of all partial derivatives of order $2$ of $\gamma$ taken, this time, with respect
   to variables $a,b,c$.
   
   Let $\Lambda$ be the linear system of curves of degree $5$ generated by
   $$u_1=a^2(5c^3-a^3),\;\; u_2=b^2(b^3-5c^3),\;\; u_3=c^2(c^3-5a^3),$$
   $$u_4=c^2(5b^3-c^3),\;\; u_5=5b^2(a^3-c^3),\;\; u_6=5a^2(c^3-b^3).$$
   These are the coefficients in \eqref{eq: gamma}.
   
   We establish the following easy fact.
\begin{lemma}\label{lem: deg 5 bpf}
   The system $\Lambda$ is base point free.
\end{lemma}   
\proof
   The polynomial $u_6$ vanishes if $a=0$ or $c=\eps^\alpha b$ for some $\alpha\in\left\{0,1,2\right\}$.
   
   Assume first that $a=0$. Then $u_3=c^5$, so it must be $c=0$. But then $u_2=b^5$ and it must be $b=0$,
   which is not possible for coordinates of a point in the projective plane.
   
   Now we turn to the second case $c=\eps^\alpha b$. Then vanishing of $u_2$ implies that $b=0$ and
   proceeding as in the previous case we conclude that all coordinates vanish.
\endproof
   Thus $\Lambda$ is not defined by vanishing along a specific set of points in $\P^2$.
   Nevertheless it has the unexpected property of having a member with a point of multiplicity $3$
   in a general point, which does not follow by a naive dimension count. It would be very interesting
   to understand better how $\Lambda$ appears and to investigate properties of the \emph{companion surface},
   i.e., the image of $\P^2$ in $\P^5$ under the morphism defined by $\Lambda$, see \cite[Section 4.4]{Szp19}
   for this path of thought.


\section{General case}
   In the present section we generalize the construction exhibited in Section \ref{sec: initial}.
   Recall, that the Fermat arrangement of lines is defined by linear factors of the polynomial
   $$F_m=(x^m-y^m)(y^m-z^m)(z^m-x^m).$$
   The singular points $S_m$ of this arrangement (i.e., points where $2$ or more lines intersect) are the union of
   a complete intersection grid $W_m$ defined by the ideal $I(W_m)=\left((x^m-y^m),(x^m-z^m)\right)$
   and the three coordinate points $X=\left\{X_1,X_2,X_3\right\}$.
   We consider sets of points $Z_m$ defined as the set difference of $S_{2m}$ and $W_m$ (therefore the name diminished Fermat-type arrangements). It turns out that these sets behave surprisingly
   regularly.
\begin{lemma}\label{lem: 3m^2 points}
   The points in $Y_m=W_{2m}\setminus W_m$ form an almost
   complete intersection. The ideal $I(Y_m)$ of $Y_m$ is generated in degree $2m$ by the following $3$ polynomials
   $$f_1=x^{2m}-y^{2m},\;\; f_2=x^{2m}-z^{2m},\;\; f_3=(x^m+z^m)(y^m+z^m).$$
\end{lemma}
\proof
   The proof is the same (with obvious exponent changes) as of Lemma \ref{lem: 27 points}
   and therefore omitted here.
\endproof
\begin{lemma}\label{lem: 3+3m^2 points}
   The ideal $I(Z_m)$ is generated in degree $2m+1$ by the following polynomials
   $$h_1=x(y^{2m}-z^{2m}),\;\; h_2=y(x^{2m}-z^{2m}),\;\; h_3=z(x^m+y^m)(y^m+z^m),$$
   $$h_4=z(x^m+y^m)(x^m+z^m),\;\; h_5=y(x^m+z^m)(y^m+z^m),\;\; h_6=x(x^m+z^m)(y^m+z^m).$$
\end{lemma}
\proof
   The proof is almost verbatim as that of Lemma \ref{lem: 30 points}.
   We omit it here.
\endproof
   The main result of this work is the following.
\begin{theorem}\label{thm: main}
   For $m\geq 3$ the sets $Z_m$ have the $U(2m+1,3)$ property.
\end{theorem}
\proof
   It is enough to produce the equation of the unexpected curve of degree $2m+1$ explicitly.
   To this end, let $P=(a:b:c)$ be a general point in $\P^2$ and let $\Gamma_m$ be the curve defined
   by
   \begin{align*}
      \gamma_m=\gamma_{m,P}(x:y:z) &= a^{m-1}((2m-1)c^m-a^m)h_1 + b^{m-1}(b^m-(2m-1)c^m)h_2\\
                      &\qquad + c^{m-1}(c^m-(2m-1)a^m)h_3 + c^{m-1}((2m-1)b^m-c^m)h_4\\
                      &\qquad + (2m-1)b^{m-1}(a^m-c^m)h_5 + (2m-1)a^{m-1}(c^m-b^m)h_6.
   \end{align*}
   Obviously $\gamma_m\in I(Z_m)$. Moreover $\mult_P\Gamma_m\geq 3$. This can be checked
   directly computing partial derivatives of $\gamma_m$ up to order $2$ and checking that
   they all vanish at $P$. By a multiple use of the Euler formula we can justify
   that claim -- we omit here simple calculations.
\endproof
\subsection{The BMSS dual curve}
   Of course, also in the general case, we can look for the equation of $\Gamma_m$
   from the perspective of coordinates $(a:b:c)$. We obtain a curve of degree $2m-1$
   \begin{align*}
      \gamma_m=\gamma_Q(a:b:c) &= x(z^{2m}-y^{2m})a^{2m-1}+ y(x^{2m}-z^{2m})b^{2m-1}+ z(y^{2m}-x^{2m})c^{2m-1} \\
                     &\qquad + (2m-1)y(x^m+z^m)(y^m+z^m)a^mb^{m-1} \\
                     &\qquad - (2m-1)x(x^m+z^m)(y^m+z^m)a^{m-1}b^m \\
                     &\qquad - (2m-1)z(x^m+y^m)(z^m+y^m)a^mc^{m-1} \\
                     &\qquad + (2m-1)x(x^m+y^m)(z^m+y^m)a^{m-1}c^m \\
                     &\qquad + (2m-1)z(x^m+y^m)(x^m+z^m)b^mc^{m-1} \\
                     &\qquad - (2m-1)y(x^m+y^m)(x^m+z^m)b^{m-1}c^m
   \end{align*}
   with a point of multiplicity $3$ at $Q=(x:y:z)$.
   
   As in the initial case, there is a full analogy with Lemma \ref{lem: deg 5 bpf}.
\begin{lemma}
   The linear system $\Lambda_m$ generated by
   $$u_1=a^{m-1}((2m-1)c^m-a^m),\;\; u_2=b^{m-1}(b^m-(2m-1)c^m),\;\; u_3=c^{m-1}(c^m-(2m-1)a^m),$$
   $$u_4=c^{m-1}((2m-1)b^m-c^m),\;\; u_5=(2m-1)b^{m-1}(a^m-c^m),\;\; u_6=(2m-1)a^{m-1}(c^m-b^m)$$
   is base point free.
\end{lemma}   
   The existence of a member of $\Lambda_m$ vanishing to order $3$ at a general point is unexpected.
   It would be interesting to explore linear systems of this kind in a more systematic way. We formulate
   it as an open problem.
\begin{problem}
   Find more examples of base points free (even very ample) linear systems admitting members
   with exceptional high multiplicity at a general point.
\end{problem}   


\paragraph*{Acknowledgement.}
   The authors thank Piotr Pokora, Tomasz Szemberg and Justyna Szpond for introducing them to this circle
   of ideas and for helpful conversations. The first author was supported by the National Science Centre 
   grant UMO 2018/31/N/ST1/02101.



\bibliographystyle{abbrv}
\bibliography{master}





\bigskip
\bigskip
\small

\bigskip
\noindent
   Jakub Kabat\\
   Department of Mathematics, Pedagogical University of Cracow,
   Podchor\c a\.zych 2,
   PL-30-084 Krak\'ow, Poland.

\nopagebreak
\noindent
   \textit{E-mail address:} \texttt{jakub.kabat@up.krakow.pl}\\

\bigskip
\noindent
   Beata Strycharz-Szemberg\\
   Department of Mathematics,
   Cracow University of Technology,
   Warszawska 24,
   PL-31-155 Krak\'ow, Poland.

\nopagebreak
\noindent
   \textit{E-mail address:} \texttt{szemberg@pk.edu.pl}\\

\end{document}